\documentclass[11pt, a4paper, reqno]{amsart}
\usepackage{amsmath}
\usepackage{geometry}
\usepackage{enumerate}
\usepackage{amsmath}
\usepackage{hyperref}
\hypersetup{colorlinks=true, linkcolor=black, citecolor=black}
\usepackage{amsfonts}
\usepackage{mathtools}
\usepackage{amssymb}
\usepackage{orcidlink}
\usepackage{tikz}
\usetikzlibrary{cd}
\usepackage{url}
\usepackage{bbding}
\usepackage{bm}
\usepackage{extarrows}
\usepackage{setspace}
\usepackage{pifont}
\usepackage{amsmath,amsthm,amssymb,amsfonts}
\usepackage{enumitem}
\numberwithin{equation}{section}
\theoremstyle{plain}
\newtheorem{theorem}{Theorem}[section]
\newtheorem{definition}[theorem]{Definition}
\newtheorem{proposition}[theorem]{Proposition}
\newtheorem{lemma}[theorem]{Lemma}
\newtheorem{question}[theorem]{Question}

\newtheorem{corollary}[theorem]{Corollary}
\newtheorem{main}{Main Theorem}

\newtheorem{remark}[theorem]{Remark}

\usepackage{amsfonts,amssymb}
\begin{document}
	\title{Normalike operators}
	\author{Youqing Ji}
	\address{School of Mathematics, Jilin University, Changchun 130012, China}
	\email{jiyq@jlu.edu.cn}
	\author{Hepu Zhang\orcidlink{0009-0009-9139-0887}}
	\address{School of Mathematics, Jilin University, Changchun 130012, China}
	\email{zhanghp25@mails.jlu.edu.cn}
	\date{}
	\thanks{Supported by the National Natural Science Foundation of China (No. 12271202)}
	
	\maketitle
	\begin{abstract}
		In this paper, for a bounded linear operator on an infinite dimensional complex Banach space, we introduce the concept of the \textit{quasi-spectrum}. At the same time, a new class of operators, called \textit{normalike} operators, is introduced as well. We primarily investigate the relationship between the quasi-spectrum and the invariant subspace problem. At first, we prove that the quasi-spectrum is a subset of $[0,1]$ containing 1 and it is right-closed. Then, we show that if the quasi-spectrum for a normalike operator is not a singleton, then the operator has a nontrivial invariant subspace. Finally, under the same conditions, we also show that the operator possesses a nontrivial hyperinvariant subspace.
	\end{abstract}

	\section{Introduction}
	Let $X$ be an infinite dimensional complex Banach space, and let $\mathcal{B}(X)$ be the algebra of all bounded linear operators on $X$. A nonempty subset $Y\subset X$ is called a linear manifold if $\alpha x+\beta y\in Y$ for all $x,y\in Y$ and $\alpha,\beta\in\mathbb C$. A closed linear manifold is called a subspace. A subspace $Y\subset X$ is called an invariant subspace for $T\in\mathcal{B}(X)$ if $Tx\in Y$ for all $x\in Y$. It is called a hyperinvariant subspace for $T$ if $Y$ is an invariant subspace for every $S\in\mathcal{B}(X)$ which commutes with $T$. A subspace is called nontrivial if it is different from both \(\{0\}\) and \(X\). Throughout this paper, $\mathbb{N}$, $\mathbb{R}$, and $\mathbb{C}$ denote the sets of positive integers, real numbers, and complex numbers, respectively.
	
	The invariant subspace problem asks whether every bounded linear operator on an infinite dimensional Banach space has a nontrivial invariant subspace. In 1975, Enflo constructed the first example of a bounded linear operator on a Banach space with no nontrivial invariant subspaces. The paper had existed in manuscript form for about 12 years prior to its publication \cite{MR892591}. In 1984, Read \cite{MR749447} produced more examples of operators on some non-reflexive Banach spaces with no nontrivial invariant subspace. These counterexamples motivate the search for conditions under which a bounded linear operator on an infinite dimensional Banach space admits a nontrivial invariant subspace.
	
	In 1954, Aronszajn and Smith \cite{MR65807} proved that there exists a nontrivial invariant subspace for any nonzero compact operator on a complex Banach space. In 1973, Lomonosov \cite{MR420305} gave a brief and elegant proof (using the Schauder fixed point theorem) of the following generalization of the above theorem: every nonzero compact operator on a complex Banach space has a nontrivial hyperinvariant subspace.

	Spectral theory and local spectral theory are also closely related to the invariant subspace problem. Every contraction operator on a complex Hilbert space whose spectrum contains the unit circle has a nontrivial invariant subspace \cite{MR2760647} (the complete proof is in \cite{MR990552}). Every local spectral subspace of an operator is a hyperinvariant linear manifold \cite{MR1747914}.

	In addition, the invariant subspace problem can be reformulated in terms of cyclic subspaces. For $x$ in $X$, let $M_x=\bigvee_{n\ge0}\left\{T^nx\right\}$, i.e., the closed linear span of $\left\{T^nx:n\ge0\right\}$. Then $M_x$ is the smallest invariant subspace containing $x$. Thus, \(T\) has a nontrivial invariant subspace if and only if \(M_x\ne X\) for some nonzero \(x\in X\). This leads to the following question:
	\begin{question}\label{que1}
		Under what conditions does there exist a nonzero vector \(x_0\in X\) such that \(M_{x_0}\ne X\)?
	\end{question}
	
	Local spectral theory provides one perspective on Question \ref{que1}. We begin by recalling some relevant notions. For $T\in\mathcal{B}(X)$, denote the spectrum and the spectral radius of $T$ by $\sigma(T)$ and $r(T)$, respectively. The spectral radius formula gives
	$$
	r(T)=\lim\limits_{n\to\infty}\Vert T^n\Vert^{1/n},
	$$
	and the local spectral radius (see \cite{MR897643}) of $T$ at $x\in X$ is defined as
	$$
	r_T(x)=\limsup\limits_{n\to\infty}\Vert T^nx\Vert^{1/n},
	$$
	where the limit superior is used because the sequence $\left\{\Vert T^nx\Vert^{1/n}\right\}_{n=1}^\infty$ may not converge. For example, let \(T\) be the backward shift operator on $\ell^2$. Define $a_0=1$ and, for $n\ge1$, set $a_n=e^{-2^k}$ whenever $2^{k-1}\le n<2^k$, $k\ge1$. Let $x=\left\{\sqrt{a^2_n-a^2_{n+1}}\right\}_{n=0}^\infty$. Then $x\in\ell^2$ and $\Vert T^nx\Vert=a_n$ for $n\ge1$. Since $$\liminf\limits_{n\to\infty}a^{1/n}_n=\mathrm{e}^{-2},\ \limsup\limits_{n\to\infty}a^{1/n}_n=\mathrm{e}^{-1},$$ $\left\{\Vert T^nx\Vert^{1/n}\right\}_{n=1}^\infty$ does not converge (see \cite{MR1264802}).
	\vspace{\baselineskip}
	
	Next, for $0\neq x\in X$, let $R^x=\sup\{r_T(y):y\in M_x\}$, and denote the restriction of $T$ to $M_x$ by $T|_{M_x}$. Applying the identity $r(T)=\sup\limits_{y\in X}r_T(y)$ (see \cite{MR1747914}) to $T|_{M_x}$, we obtain $R^x=r(T|_{M_x})$. If $r(T|_{M_x})=r_T(x)$ for every nonzero $x\in X$, and if there exist $x,\ x'\in X\setminus\{0\}$ such that $r_T(x)\neq r_T(x')$, then at least one of \(M_x\) and \(M_{x'}\) is a nontrivial invariant subspace for \(T\). Note that $\left\{r_T(x):0\neq x\in X\right\}$ may be a singleton set. For example, $\left\{r_T(x):0\neq x\in X\right\}=\{0\}$ for a quasinilpotent operator $T$ (i.e., $r(T)=0$). Consequently, two questions arise: 
	\begin{enumerate}[label=(\arabic*)]
		\item Under what conditions does \(\{r_T(x):0\ne x\in X\}\) contain more than one point?
		\item Under what conditions does \(r(T|_{M_x})=r_T(x)\) hold for every nonzero \(x\in X\)?
	\end{enumerate}
	
	\vspace{\baselineskip}
	
In 1995, Atzmon \cite{MR1264802} introduced the concept of \textit{power regular} operators. Atzmon's results provide a partial answer to the second question above. An operator $T\in\mathcal{B}(X)$ is called \textit{power regular} if, for every $x\in X$, the limit $\lim\limits_{n\to\infty}\Vert T^nx\Vert^{1/n}$ exists. Isometric operators and quasinilpotent operators are power regular. Atzmon established a criterion ensuring that \(T\) is power regular and that $ r_T(x)=r(T|_{M_x})$ for every $0\ne x\in X$. Atzmon also showed that this criterion is satisfied by every decomposable operator and by every operator \(T\) for which the set $\{|\lambda|:\lambda\in\sigma(T)\}$ has empty interior in \([0,\infty)\). For details on decomposable operators, the reader is referred to \cite{MR152893}. Later, Bourdon \cite{MR1460419} proved that every hyponormal operator $T$ on a complex Hilbert space $\mathcal{H}$ is power regular, and that $T$ has a nontrivial invariant subspace if there exists $0\ne x\in\mathcal{H}$ such that
	$r_T(x)<\Vert T\Vert$. In 2025, Hu and Ji \cite{MR4905072} proved that the inverse of an invertible power regular operator is not necessarily power regular.

	Nevertheless, the preceding argument based on distinct local spectral radii yields no information when \(T\) is quasinilpotent, since \(r_T(x)=0\) for every \(x\in X\).
	
	We next present two approaches to the invariant subspace problem for quasinilpotent operators. One approach, based on rank-one perturbations, was developed by Tcaciuc in 2019 (see \cite{MR3950035}). Tcaciuc proved that a quasinilpotent operator $T$ has a nontrivial invariant subspace if and only if there exist a rank-one operator $F$ and $\alpha\in\mathbb{C}\setminus\{0,1\}$ such that $T+F$ and $T+\alpha F$ are also quasinilpotent. Under additional assumptions, Tcaciuc further proved that, for every rank-one operator $F$, the perturbation $T+\alpha F$ has a invariant subspace of infinite dimension and codimension for all nonzero $\alpha\in\mathbb C$, with at most two exceptions.
	
	Another approach involves the power set of a quasinilpotent	operator, introduced by Douglas and Yang \cite{MR3837636} in 2018. This provides a tool for studying hyperinvariant subspaces of such operators.
	
	Let $T\in\mathcal{B}(X)$ be a quasinilpotent operator. Define
	$$
	k_T(x):=\limsup_{z\to 0}\dfrac{\ln\Vert(T-z)^{-1}x\Vert}{\ln\Vert(T-z)^{-1}\Vert},\quad 0\neq x\in X,
	$$
	and call $\Lambda(T):=\left\{k_T(x):0\neq x\in X\right\}$ the power set of $T$. Then $\Lambda(T)\subset[0,1]$. For $\tau\in[0,1]$, let
	$$
	M_\tau:=\left\{0\neq x\in X:k_T(x)\le\tau\right\}\cup\{0\}.
	$$
	
	The following result was established in \cite{MR4211134}.
	
	\begin{proposition}
		Let $T\in\mathcal{B}(X)$ be a quasinilpotent operator, and let $\tau\in[0,1]$. If $M_\tau$ is closed, it is a hyperinvariant subspace for $T$. In particular, if there exist $\tau_1\neq\tau_2\in\Lambda(T)$ such that both $M_{\tau_1}$ and $M_{\tau_2}$ are closed, then $T$ has a nontrivial hyperinvariant subspace.
	\end{proposition}
	
	It is natural to ask whether $M_\tau$ is always closed. In \cite{MR4260659}, Ji and Liu constructed a quasinilpotent operator $T$ with $\Lambda(T)=[0,1]$ such that, for every $\tau\in[0,1)$, $M_\tau$ is not closed.
	
	Following the idea of focusing on whether $\left\{r_T(x):0\neq x\in X\right\}$ or $\Lambda(T)$ is a singleton, we introduce a concept called the \textit{quasi-spectrum}. And we are still interested in what happens when the quasi-spectrum is not a singleton. In Section \ref{sec3}, we will give examples of quasinilpotent operators whose quasi-spectra are not singletons.
	
	\begin{definition}
		Given $T\in\mathcal{B}(X)$ be non-nilpotent, and a sequence of positive numbers $\{r_n\}_{n\in\mathbb{N}}$ with $\lim\limits_{n\to\infty}\Vert r_nT^n\Vert^{1/n}=1$, define $$q_T(x)=\limsup\limits_{n\to\infty}\Vert r_nT^nx\Vert^{1/n},\ \forall x\in X,$$ and denote by $\sigma_q(T)$ the set $\left\{q_T(x):0\neq x\in X\right\}$, which is called the quasi-spectrum for $T$.
	\end{definition}
	
	One may notice the symbol we use does not contain $\{r_n\}$. This is bacause the value of $q_T(x)$ is independent of the choise of $\{r_n\}_{n\in\mathbb{N}}$. In fact, supposing that$$\lim\limits_{n\to\infty}\Vert r_nT^n\Vert^{1/n}=\lim\limits_{n\to\infty}\Vert r'_nT^n\Vert^{1/n}=1,$$ then $\lim\limits_{n\to\infty}\left(\dfrac{r_n}{r'_n}\right)^{1/n}=1$. Hence, $$\limsup\limits_{n\to\infty}\Vert r_nT^nx\Vert^{1/n}=\limsup\limits_{n\to\infty}\left(\Vert r'_nT^nx\Vert^{1/n}\cdot\left(\dfrac{r_n}{r'_n}\right)^{1/n}\right)=\limsup\limits_{n\to\infty}\Vert r'_nT^nx\Vert^{1/n}.$$Therefore, setting \(r_n=\|T^n\|^{-1}\) may simplify the computations.
	
	In addition, it is clear that $q_T(x)$ is independent of $\|x\|$, i.e., $q_T(ax)=q_T(x)$, for every $a\in\mathbb{C}\setminus\{0\}$ and $x\in X$. Thus, $\sigma_q(T)=\left\{q_T(x):x\in X,\ \Vert x\Vert=1\right\}$.
	
	Besides, for the local spectral radius, $0\le r_T(x)\le r(T)$ for all $x\in X$ and there exists an $x_0\in X$ such that $r_T(x_0)=r(T)$ (see \cite{MR897643}). And it is the same for power set; that is, $0\le k_T(x)\le1$ for all $x\in X$ and there is an $x_0\in X$ such that $k_T(x_0)=1$ (see \cite{hu2026powersetquasinilpotentoperators}). Similarly, we shall prove that $0\le q_T(x)\le1$ for all $x\in X$ and there must be an $x_0\in X$ such that $q_T(x_0)=1$.
	
	\begin{definition}\cite{MR4654136}
		A nonempty set $E\subset\mathbb{R}$ is called right-closed if $\sup F\in E$ for every nonempty bounded subset $F$ of $E$.
	\end{definition}
	
	It was proved that the power set is right-closed in \cite{MR4654136}. We will show that the quasi-spectrum is also right-closed in Section \ref{sec2}.
	
	The following question is of interest to us.
	
	\begin{question}\label{que}
		When the quasi-spectrum of a given operator is not a singleton, does it have a nontrivial invariant subspace? Moreover, does it have a nontrivial hyperinvariant subspace?
	\end{question}
	
	To answer Question \ref{que}, we shall introduce a new class of operators (on Banach spaces), which is a generalization of normal operators (on Hilbert spaces). Several generalizations of normal operators have been studied, including subnormal operators and hyponormal operators. A bounded linear operator $T$ on a complex Hilbert space $\mathcal{H}$ is said to be subnormal if there exists a normal operator $N$ on a Hilbert space $\mathcal{K}\supseteq\mathcal{H}$ such that $\mathcal{H}$ is invariant under $N$ and $T=N|_{\mathcal{H}}$. Moreover, $T$ is said to be hyponormal if $T^\ast T-TT^\ast\ge 0$. In 1978, Brown \cite{MR511974} proved that each subnormal operator has a nontrivial invariant subspace. Later in 1987, Brown \cite{MR873378} proved that if $T$ is hyponormal and $R(\sigma(T))\neq C(\sigma(T))$, then $T$ has a nontrivial invariant subspace, where $C(\sigma(T))$ denotes the space of all continuous functions on $\sigma(T)$, and $R(\sigma(T))$ denotes the closure in $C(\sigma(T))$ of the rational functions with poles off $\sigma(T)$.
	
	Before introducing our new generalization, we briefly recall a related result concerning normal operators. By the spectral theorem, every normal operator is unitarily equivalent to a multiplication operator on an $L^2$ space. More precisely, for a normal operator $N$ on a Hilbert space $\mathcal{H}$, there exist a measure space $(\Omega,\mu)$, a function $g\in L^\infty(\Omega,\mu)$, and a unitary operator
	$U:\mathcal{H}\to L^2(\Omega,\mu)$ such that $UNU^{-1}=T_g$, where $T_gf=gf$ for $f\in L^2(\Omega,\mu)$. Fix $x\in\mathcal{H}$ with $\Vert x\Vert=1$. Set $f=Ux$. Then we have $\Vert N^nx\Vert^{1/n}=\Vert T_g^nf\Vert^{1/n} (\forall n\in\mathbb{N})$ and $\Vert f\Vert=1$. Define a measure $\nu$ such that its Radon-Nikodym derivative with respect to $\mu$ is $|f|^2$, i.e., $\mathrm{d}\nu=|f|^2\,\mathrm{d}\mu$, and hence $\nu$ is a probability measure. Given $n\in\mathbb{N}$, let $p=\dfrac{n+1}{n}$, and $p'=n+1$. By H\"older inequality, we obtain$$\begin{aligned}
	\Vert N^nx\Vert^{1/n}&=\left(\int_{\Omega}|g(t)|^{2n}|f(t)|^2\,\mathrm{d}\mu\right)^{1/2n}\\
	&=\left(\int_{\Omega}|g(t)|^{2n}\,\mathrm{d}\nu\right)^{1/2n}\\
	&\le\left(\int_{\Omega}|g(t)|^{2np}\,\mathrm{d}\nu\right)^{1/2np}\cdot\left(\int_{\Omega}1^{p'}\,\mathrm{d}\nu\right)^{1/2np'}\\
	&=\left(\int_{\Omega}|g(t)|^{2(n+1)}\,\mathrm{d}\nu\right)^{1/2(n+1)}\\
	&=\Vert N^{n+1}x\Vert^{1/(n+1)}.
	\end{aligned}$$ Therefore, $\left\{\Vert N^nx\Vert^{1/n}\right\}_{n=1}^\infty$ is increasing in $n$ for every unit vector $x\in\mathcal{H}$. Furthermore, one can obtain
	\begin{equation}\label{gen}
		r_N(x)=\limsup_{n\to\infty}\Vert N^nx\Vert^{1/n}=\sup_{n\in\mathbb{N}}\left\{\Vert N^nx\Vert^{1/n}\right\}, \tag{$\ast$}
	\end{equation}
	for every unit vector $x\in\mathcal{H}$, and our generalization is based on \eqref{gen}. 
	\vspace{\baselineskip}
	
	Moreover, since $\left\{\Vert N^nx\Vert^{1/n}\right\}_{n=1}^\infty$ is increasing, $\Vert Nx\Vert\le\Vert N^nx\Vert^{1/n}$ for all $n$. Then $\Vert Nx\Vert\le r_N(x)$ for every unit vector $x\in\mathcal{H}$. For $y\in\operatorname{span}_{m\ge0}\{N^mx\}$, clearly $r_N(y)\le r_N(x)$. Hence, $$\dfrac{\Vert Ny\Vert}{\Vert y\Vert}\le r_N(y/\|y\|)=r_N(y)\le r_N(x).$$By continuity, $\Vert Ny\Vert\le r_N(x)\|y\|$ holds for every $y\in M_x$. It follows that $\|N|_{M_x}\|\le r_N(x)$. Note that $N|_{M_x}$ is subnormal. So $r(N|_{M_x})=\|N|_{M_x}\|\le r_N(x)$. Obviously, $r(N|_{M_x})\ge r_N(x)$. Therefore, $r(N|_{M_x})=r_N(x)$ holds for every nonzero $x\in\mathcal{H}$.
		
		\vspace{\baselineskip}
		
	Motivated by the idea of generalizing normal operators, we, again, introduce a new notion termed \textit{normalike} operators in this paper, which is helpful for answering Question \ref{que}.

	\begin{definition}\label{like}
		Let $T\in\mathcal{B}(X)$ be non-nilpotent. Say $T$ is normalike if there exists a sequence of positive numbers $\{r_n\}_{n\in\mathbb{N}}$ with $\lim\limits_{n\to\infty}\Vert r_nT^n\Vert^{1/n}=1$ such that $$q_T(x)=\sup_{n\in\mathbb{N}}\left\{\Vert r_nT^nx\Vert^{1/n}\right\}$$
		for all unit vectors $x$ in $X$.
	\end{definition}

	Following the previous discussion, for a nonzero normal operator $N$ on a Hilbert space $\mathcal{H}$, $r_N(x)=\sup_{n\in\mathbb{N}}\left\{\Vert N^nx\Vert^{1/n}\right\}$ for every unit vector $x\in\mathcal{H}$. Taking $r_n=\Vert N\Vert^{-n}$, we conclude that $N$ is normalike. Actually, every nonzero hyponormal operator is also normalike, since the equality $r_T(x)=\sup_{n\in\mathbb{N}}\left\{\Vert T^nx\Vert^{1/n}\right\}$ also holds for a hyponormal operator $T$ (one can see \cite{MR1460419}). In Section \ref{sec3}, we will give examples of normalike operators which are not hyponormal.

	Our main theorems are as follows. And Main theorem \ref{inv} and Main theorem \ref{hyp} partially answer Question \ref{que1} and Question \ref{que}.
	
	\begin{main}\label{1111}
		Let $T\in\mathcal{B}(X)$ be non-nilpotent. Then the followings hold:\begin{enumerate}[label=(\arabic*)]
			\item $\{1\}\subset \sigma_q(T)\subset[0,1]$;
			\item $\sigma_q(T)$ is right-closed.
		\end{enumerate}
	\end{main}
	
	\begin{main}\label{inv}
		Suppose that $T\in\mathcal{B}(X)$ is normalike. If there exists an $x_0\in X$ such that $q_T(x_0)<1$, then $M_{x_0}$ is a nontrivial invariant subspace. Hence, if $\sigma_q(T)\neq\{1\}$, then $T$ has a nontrivial invariant subspace.
	\end{main}

	\begin{main}\label{hyp}
		Suppose that $T\in\mathcal{B}(X)$ is normalike. If $\sigma_q(T)\neq\{1\}$, then $T$ has a nontrivial hyperinvariant subspace.
	\end{main}

	\section{The quasi-spectrum}\label{sec2}
	
	In this section, we prove Main Theorem \ref{1111}: the quasi-spectrum is a subset of $[0,1]$, contains 1 and is right-closed. To prove $1\in\sigma_q(T)$, we need the following lemma.

	\begin{lemma}\cite{MR1847174}\label{1inQ}
		Suppose $T\in\mathcal{B}(X)$ and given a sequence of positive numbers $\{a_n\}_{n=1}^\infty$ with $a_n\to0$. For $x\in X$, let $\Omega_x=\left\{n\in\mathbb{N}:\Vert T^nx\Vert\ge a_n\Vert T^n\Vert\right\}$, and let $\Omega=\left\{x\in X:\Omega_x\text{ is infinite}\right\}$. Then $\Omega$ is of the second category.
	\end{lemma}
	
	\begin{proof}[\textbf{Proof of Main Theorem \ref{1111} (1)}]
	\begin{enumerate}[label=\arabic*)]
		\item Fix $x\in X\setminus\{0\}$. Take a sequence of positive numbers $\{r_n\}_{n\in\mathbb{N}}$ with $\lim\limits_{n\to\infty}\Vert r_nT^n\Vert^{1/n}=1$. Then $$0\le q_T(x)=\limsup\limits_{n\to\infty}\Vert r_nT^nx\Vert^{1/n}\le\limsup\limits_{n\to\infty}\Vert r_nT^n\Vert^{1/n}\Vert x\Vert^{1/n}=1.$$ So $\sigma_q(T)\subset[0,1]$.
		
		\item Since the value of $q_T(x)$ is independent of the choise of $\{r_n\}_{n\in\mathbb{N}}$, we take $r_n=\|T^n\|^{-1}\ (n\in\mathbb{N})$. Let $a_n=n^{-1}\ (n\in\mathbb{N})$, and thus $a_n\to0(n\to\infty)$. According to Lemma \ref{1inQ}, there exists some $x_0\in X\setminus\{0\}$ such that $\Vert r_nT^nx_0\Vert\ge r_na_n\Vert T^n\Vert=n^{-1}$ holds for infinitely many $n$. Then there exist a subquence $\{n_k\}_{k=1}^\infty$ such that $\Vert r_{n_k}T^{n_k}x_0\Vert\ge n^{-1}_k$ for all $k\in\mathbb{N}$. Since $n^{1/n}\to1(n\to\infty)$, $\Vert r_{n_k}T^{n_k}x_0\Vert^{1/n_k}\to1(k\to\infty)$. It follows that $q_T(x_0)=1$. Therefore, $1\in \sigma_q(T)$.
			\end{enumerate}
	\end{proof}
	
	Next we show that $\sigma_q(T)$ is right-closed. Recall that the meaning of $\sigma_q(T)$ is right-closed is that $\sup F\in\sigma_q(T)$ for every nonempty subset $F$ of $\sigma_q(T)$. If $\sigma_q(T)$ is finite, it is obviously right-closed. Next, assume that $\sigma_q(T)$ is infinite, and let $F\subset\sigma_q(T)$ be nonempty. If $\sup F\in F$, then we are done. If $\sup F\notin F$, then the approach to proving $\sup F\in E$ is to take an infinite strictly increasing sequence in $F$ which converges to $\sup F$ as an auxiliary tool. The detailed proof is as follows.

		\begin{proof}[\textbf{Proof of Main Theorem \ref{1111} (2)}]
		Suppose that $\sigma_q(T)$ is infinite. Based on the above discussion, in order to prove $\sigma_q(T)$ is right-closed, it suffices to prove that $t\in\sigma_q(T)$ whenever $\{t_k\}_{k=1}^\infty\subset\sigma_q(T)$ is a strictly increasing sequence with the limit $t\in(0,1]$. Take a sequence of positive numbers $\{r_n\}_{n\in\mathbb{N}}$ with $\lim\limits_{n\to\infty}\Vert r_nT^n\Vert^{1/n}=1$ and take $\{x_k\}_{k=1}^\infty\subset X$ such that $q_T(x_k)=t_k$ for every $k\in\mathbb{N}$. When $t=1$, the conclusion follows directly from part (1) of Main Theorem \ref{1111}. Hence, we assume that $t\in(0,1)$, and we prove $t\in\sigma_q(T)$ step by step.
		\\
		\\
		\noindent\textbf{Step 1:}
		Let $r_0=1$. Define
		\[
		Y:=
		\left\{
		x\in X:
		\sup_{n\ge0}t^{-n}\|r_nT^n x\|<+\infty
		\right\},
		\qquad
		\|x\|_Y:=
		\sup_{n\ge0}t^{-n}\|r_nT^n x\|.
		\]
		\\
		\\
			\noindent\textbf{Step 2:} Claim that $(Y,\|\cdot\|_Y)$ is a Banach space.
		\\
		\\
	\noindent	\textbf{Step 2.1: $Y$ is a linear space.}
	
	At first, clearly $0\in Y$, and thus $Y$ is not empty. Let $x,\ y\in Y$, and let $\alpha\in\mathbb{C}$. Then we have $$\begin{aligned}
		\sup_{n\ge0}t^{-n}\|r_nT^n(x+y)\|&\le\sup_{n\ge0}\left(t^{-n}\|r_nT^n x\|+t^{-n}\|r_nT^n y\|\right)\\
		&\le\sup_{n\ge0}t^{-n}\|r_nT^n x\|+\sup_{n\ge0}t^{-n}\|r_nT^n y\|<+\infty.
	\end{aligned}$$
	It implies that $x+y\in Y$. Since $$\sup_{n\ge0}t^{-n}\|r_nT^n(\alpha x)\|=|\alpha|\cdot\sup_{n\ge0}t^{-n}\|r_nT^n x\|<+\infty,$$ it follows that $\alpha x\in Y$. So $Y$ is a linear space.
		\\
		\\
		\noindent	\textbf{Step 2.2: $\|\cdot\|_Y$ is a norm which is stronger than $\|\cdot\|$.}
		
		Since the supremum includes $n=0$, $\|x\|_Y\ge\|x\|$ for every $x\in Y$. This implies positive definiteness of $\|\cdot\|_Y$. By Step 2.1, $\|\cdot\|_Y$ is homogeneous and it also satisfies the triangle inequality. Hence $\|\cdot\|_Y$ is a norm and it is stronger than $\|\cdot\|$.
		\\
		\\
		\noindent	\textbf{Step 2.3: $Y$ is complete.}
		
		Let $\{y_j\}_{j\ge1}\subset Y$ be a Cauchy sequence in
		$\|\cdot\|_Y$. Since $\|y_j-y_\ell\|\le\|y_j-y_\ell\|_Y$, $\{y_j\}_{j\ge1}$ is also a Cauchy sequence in $X$.
		Since $X$ is complete, there exists $y\in X$ such that $\|y_j-y\|\to0 (j\to\infty)$.
		Then we need to establish that $y\in Y$ and $\|y_j-y\|_Y\to0 (j\to\infty)$.
		
		Fix $\varepsilon>0$. Since $\{y_j\}$ is Cauchy in $Y$, there
		exists $J_\varepsilon\ge1$ such that $\|y_j-y_\ell\|_Y\le\varepsilon$ whenever $j,\ell\ge J_\varepsilon$.
		Equivalently,
		$$t^{-n}\|r_nT^n(y_j-y_\ell)\|\le\varepsilon
			\quad
			\text{for every }n\ge0
			\text{ and all }j,\ell\ge J_\varepsilon.$$
		
		Fix $j\ge J_\varepsilon$ and $n\ge0$.
		Since $T^n$ is a bounded linear operator on $X$,
	$$\lim\limits_{\ell\to\infty}t^{-n}\|r_nT^n(y_j-y_\ell)\|=t^{-n}\|r_nT^n(y_j-y)\|.$$
		It follows that
		$$\sup_{n\ge0}t^{-n}\|r_nT^n(y_j-y)\|
			\le\varepsilon
			\quad\text{for every }j\ge J_\varepsilon,$$
		and hence \begin{equation}
			\label{eq:limit-bound}\lim\limits_{j\to\infty}\sup_{n\ge0}t^{-n}\|r_nT^n(y_j-y)\|=0.
		\end{equation}
		
		Since $$\sup_{n\ge0}t^{-n}\|r_nT^ny\|\le\sup_{n\ge0}t^{-n}\|r_nT^n(y_j-y)\|+\sup_{n\ge0}t^{-n}\|r_nT^ny_j\|$$ holds for every $j\in\mathbb{N}$, we obtain $$\sup\limits_{n\ge0}t^{-n}\|r_nT^ny\|<+\infty.$$ Then $y\in Y$. By \eqref{eq:limit-bound}, $\lim\limits_{j\to\infty}\|y_j-y\|_Y=0$. Therefore, $Y$ is complete, and thus $Y$ is a Banach space.
		\\
		\\
			\noindent\textbf{Step 3:} For $k,m\in\mathbb{N}$, define
		\[
		E_{k,m}:=
		\left\{
		x\in X:
		\|r_nT^n x\|\le m t_k^n
		\text{ for every }n\ge0
		\right\}.
		\]
		
		If $x\in E_{k,m}$, then $$t^{-n}\|r_nT^n x\|\le m \left(\dfrac{t_k}{t}\right)^n,$$ and thus $\sup\limits_{n\ge0}t^{-n}\|r_nT^n x\|\le m<+\infty$. It implies that $x\in Y$, and therefore, $E_{k,m}\subset Y$.
		\\
		\\
			\noindent\textbf{Step 4: Every $E_{k,m}$ is closed in $\|\cdot\|_Y$.}
		
		Select a pairwise disjoint convergent sequence $\{z_j\}_{j=1}^\infty$ in $E_{k,m}$, and let $z_0=\lim\limits_{j\to\infty}z_j$. Notice that the convergence here means $\|z_j-z_0\|_Y\to0(j\to\infty)$, and $0<t<1$. Therefore, for every $n\ge0$ and $j\in\mathbb{N}$,
		$$\begin{aligned}
			\|r_nT^nz_0\|&\le\|r_nT^nz_j\|+\|r_nT^n(z_0-z_j)\|\\
			&\le m t_k^n+t^{-n}\|r_nT^n(z_0-z_j)\|\\
			&\le m t_k^n+\sup_{n\ge0}t^{-n}\|r_nT^n(z_0-z_j)\|\\
			&=m t_k^n+\|z_j-z_0\|_Y.
		\end{aligned}$$
		
		Taking the limit as $j\to\infty$ gives $\|r_nT^nz_0\|\le m t_k^n$. It follows that $z_0\in E_{k,m}$. Hence, $E_{k,m}$ is closed in $\|\cdot\|_Y$.
		\\
		\\
		\noindent\textbf{Step 5:} $\bm{\left\{x\in X:q_T(x)<t\right\}=\displaystyle\bigcup_{k,m\in\mathbb{N}}E_{k,m}}$.
		
		On the one hand, given an $x\in X$ with $q_T(x)<t$, there is $k_0\in\mathbb{N}$ such that $q_T(x)\le t_{k_0}$. By the definition of the limit superior, for $\varepsilon=\dfrac{1}{2}\left(t_{k_0+1}-t_{k_0}\right)>0$, there exists $N=N_\varepsilon\in\mathbb{N}$ such that for every $n>N$, we have
		$$\Vert r_nT^nx\Vert^{1/n}\le t_{k_0}+\varepsilon<t_{k_0+1}.$$
		Hence, $\Vert r_nT^nx\Vert\le t^n_{k_0+1},\ \forall n>N$. For $n=1,2\cdots,N$, set $$\alpha_n=\dfrac{\Vert r_nT^nx\Vert}{t^n_{k_0+1}},$$and set $$m_n=1+\lfloor \alpha_n\rfloor,$$ where $\lfloor a\rfloor$ denotes the largest integer less than or equal to $a$. Let $m=1+\max\limits_{1\le n\le N}m_n$. Then $$\Vert r_nT^nx\Vert\le m_nt^n_{k_0+1}\le mt^n_{k_0+1},\ \forall0\le n\le N,$$ and
		$$\Vert r_nT^nx\Vert\le t^n_{k_0+1}\le mt^n_{k_0+1},\ \forall n>N.$$
		It follows that $\Vert r_nT^nx\Vert\le mt^n_{k_0+1}$ holds for all $n\ge0$, which implies $x\in E_{m,k_0+1}\subset\displaystyle\bigcup_{k,m\in\mathbb{N}}E_{k,m}$. So $$\left\{x\in X:q_T(x)<t\right\}\subset\displaystyle\bigcup_{k,m\in\mathbb{N}}E_{k,m}.$$
		
		On the other hand, given an $x\in\displaystyle\bigcup_{k,m\in\mathbb{N}}E_{k,m}$, there exist $k_0\in\mathbb{N}$ and $m_0\in\mathbb{N}$ such that $x\in E_{k_0,m_0}$. Then for every $n\ge0$, $\|r_nT^n x\|\le m_0 t_{k_0}^n$, and thus $q_T(x)\le t_{k_0}<t$; that is, $x\in\left\{x\in X:q_T(x)<t\right\}$. So $$\displaystyle\bigcup_{k,m\in\mathbb{N}}E_{k,m}\subset\left\{x\in X:q_T(x)<t\right\}.$$
		
		Therefore, $\left\{x\in X:q_T(x)<t\right\}=\displaystyle\bigcup_{k,m\in\mathbb{N}}E_{k,m}\subset Y$.
		\\
		\\
		\noindent\textbf{Step 6:} Since $\left\{x\in X:q_T(x)<t\right\}\subset Y$, $x_k\in Y$ for all $k$. The aim is to prove that there exists an $x_0\in Y$ such that $q_T(x_0)=t$. Assume, for a contradiction, that $q_T(x)<t$ for every $x\in Y$. Then $Y=\left\{x\in X:q_T(x)<t\right\}$. It follows from Step 5 that $Y=\displaystyle\bigcup_{k,m\in\mathbb{N}}E_{k,m}$. By the Baire category theorem, some $E_{k_0,m_0}$ has an interior point in $Y$; that is, there are some $a\in E_{k_0,m_0}$ and $\delta>0$ such that $B_Y(a,\delta)\subset E_{k_0,m_0}$. It implies that
		 \begin{equation}\label{6}
		 	\|r_nT^n h\|\le\|r_nT^n a\|+\|r_nT^n (a+h)\|\le 2m_0 t_{k_0}^n
		 \end{equation} 
		holds for all $n\ge0$ and every $h\in Y$ with $\|h\|_Y<\delta$. Now, for every $x\in Y$, since $$x=\dfrac{2\|x\|_Y}{\delta}\cdot\dfrac{\delta x}{2\|x\|_Y}$$ and $\left\|\dfrac{\delta x}{2\|x\|_Y}\right\|_Y=\dfrac{\delta}{2}<\delta$, it follows from inequality \eqref{6} that
		\[ \|r_nT^n x\|=\dfrac{2\|x\|_Y}{\delta}\cdot \left\|r_nT^n \left(\dfrac{\delta x}{2\|x\|_Y}\right)\right\|\le \frac{4m_0}{\delta}\|x\|_Y\,t_{k_0}^n ,\ n\ge0. \]
		Consequently, $q_T(x)\le t_{k_0}$ for every $x\in Y$. This contradicts $x_k\in Y$ and $t_k>t_{k_0}$ for every $k>k_0$. Thus there exists $x_0\in Y\subset X$ such that
		$q_T(x_0)=t=\lim\limits_{k\to\infty}t_k$. Therefore, $\sigma_q(T)$ is right-closed, and the proof is complete.
		
	\end{proof}

	\section{Examples of quasinilpotent normalike operators}\label{sec3}
	In Introduction, we mentioned every hyponormal operator is normalike. In addition to hyponormal operators, some special weighted shift operators are also normalike. For $1\le p<+\infty$, We abbreviate $\ell^p(\mathbb{N})$ as $\ell^p$. Let $\{e_k\}_{k=0}^\infty$ be the canonical basis of $\ell^p$.

	\begin{definition}\cite{MR361899}
		Let $\{\omega_k\}_{k=0}^{\infty}\subset\mathbb{C}$ be bounded and $1\le p<+\infty$. An operator $T\in\mathcal{B}(\ell^p)$ is called a forward weighted shift operator on $\ell^p$ if $Te_k=\omega_{k}e_{k+1},\ k=0,1,2,\cdots$. The sequence $\{\omega_k\}_{k=0}^{\infty}$ is called the weight sequence of $T$.
	\end{definition}
	
	We will calculate the quasi-spectra for some special forward weighted shift operators. These examples illustrates several facts:\begin{enumerate}[label=(\Roman*)]
		\item The class of hyponormal operators is a proper subset of the class of normalike operators.\label{sub}
		\item The quasi-spectrum for a quasinilpotent operator may contain more than one point.\label{a-k}
		\item In the Hilbert space setting, while $\sigma(T^\ast)=\left\{\overline{\lambda}:\lambda\in\sigma(T)\right\}$ holds, it fails for the quasi-spectrum.\label{spe}
		\item Every nonempty right-closed subset of $[0,1]$ containing 1 is the quasi-spectrum for some quasinilpotent operator.\label{re}
	\end{enumerate}
	
	Before investigating those, it is neccessary to make two preparations.
	
	First, for convenience of computation, the following proposition is useful. It states that the quasi-spectrum remains invariant under similarity.
	
	\begin{proposition}\label{W}
		Let $X_1$ and $X_2$ be Banach spaces and let $T_i \in \mathcal{B}(X_i)$ be non-nilpotent, $i=1,2$. If $T_1$ and $T_2$ are similar, i.e., there exists an invertible operator $W\in\mathcal{B}(X_1,X_2)$ such that $WT_1=T_2W$, then $\sigma_q(T_1)=\sigma_q(T_2)$.
	\end{proposition}
	
	\begin{proof}
		Take a sequence of positive numbers $\{r_n\}_{n\in\mathbb{N}}$ with $\lim\limits_{n\to\infty}\Vert r_nT^n_1\Vert^{1/n}=1$. Since $$\Vert r_nT^n_2\Vert^{1/n}=\Vert r_nWT^n_1W^{-1}\Vert^{1/n}\le\Vert r_nT^n_1\Vert^{1/n}\|W\|^{1/n}\|W^{-1}\|^{1/n},$$and$$\Vert r_nT^n_1\Vert^{1/n}=\Vert r_nW^{-1}T^n_2W\Vert^{1/n}\le\Vert r_nT^n_2\Vert^{1/n}\|W\|^{1/n}\|W^{-1}\|^{1/n},$$
		we have $\lim\limits_{n\to\infty}\Vert r_nT^n_2\Vert^{1/n}=1$.

		Fix $x\in X_1$. It can be calculated that
		$$\begin{aligned}
			q_{T_2}(Wx)&=\limsup_{n\to\infty}\Vert r_nT^n_2Wx\Vert^{1/n}\\
			&=\limsup_{n\to\infty}\Vert r_nWT^n_1x\Vert^{1/n}\\
			&\le\limsup_{n\to\infty}\Vert r_nT^n_1x\Vert^{1/n}\|W\|^{1/n}\\
			&=q_{T_1}(x),
		\end{aligned}$$
		and 
		$$\begin{aligned}
			q_{T_1}(x)&=\limsup_{n\to\infty}\Vert r_nT^n_1x\Vert^{1/n}\\
			&=\limsup_{n\to\infty}\Vert r_nW^{-1}WT^n_1x\Vert^{1/n}\\
			&=\limsup_{n\to\infty}\Vert r_nW^{-1}T^n_2Wx\Vert^{1/n}\\
			&\le\limsup_{n\to\infty}\Vert r_nT^n_2Wx\Vert^{1/n}\|W^{-1}\|^{1/n}\\
			&=q_{T_2}(Wx),
		\end{aligned}$$
		so we have $q_{T_1}(x)=q_{T_2}(Wx)$. Therefore, $\sigma_q(T_1)\subset\sigma_q(T_2)$. By symmetry, $\sigma_q(T_1)\supset\sigma_q(T_2)$. Hence, $\sigma_q(T_1)=\sigma_q(T_2)$.
	\end{proof}

	Since the weighted shift operator with weight sequence $\{\omega_k\}_{k=0}^{\infty}$ is similar to the weighted shift operator with weight sequence $\{|\omega_k|\}_{k=0}^{\infty}$ (see \cite{MR4493762}), we shall assume henceforth that $\omega_k\ge0(\forall k\ge0)$ by Proposition \ref{W}.
	
	Second, there is a way to determine whether an operator $T$ is normalike.
	
	\begin{lemma}
		Let $T\in\mathcal{B}(X)$ be non-nilpotent, let $c\ge0$, and let $\{r_n\}_{n\in\mathbb{N}}$ be a sequence of positive numbers with $\lim\limits_{n\to\infty}\Vert r_nT^n\Vert^{1/n}=1$. If $q_T(x_0)=c$  for some $x_0\in X$ and $\Vert r_nT^nx_0\Vert^{1/n}\le c$ for every $n\in\mathbb{N}$, then $q_T(x_0)=\sup\limits_{n\in\mathbb{N}}\left\{\Vert r_nT^nx_0\Vert^{1/n}\right\}$.	
	\end{lemma}
	
	\begin{proof}
		Since $\Vert r_nT^nx_0\Vert^{1/n}\le\sup\limits_{n\in\mathbb{N}}\left\{\Vert r_nT^nx_0\Vert^{1/n}\right\}$ holds for every $n\in\mathbb{N}$, $$q_T(x_0)\le\sup_{n\in\mathbb{N}}\left\{\Vert r_nT^nx_0\Vert^{1/n}\right\}.$$ It follows that $$c=q_T(x_0)\le\sup_{n\in\mathbb{N}}\left\{\Vert r_nT^nx_0\Vert^{1/n}\right\}\le c.$$ Thus, $q_T(x_0)=\sup_{n\in\mathbb{N}}\left\{\Vert r_nT^nx_0\Vert^{1/n}\right\}$.
	\end{proof}
	
	\begin{corollary}\label{nnnn}
		Let $T\in\mathcal{B}(X)$ be non-nilpotent and let $\{r_n\}_{n\in\mathbb{N}}$ be a sequence of positive numbers with $\lim\limits_{n\to\infty}\Vert r_nT^n\Vert^{1/n}=1$. If for every unit vector $x\in X$, there exists a positive number $c_x$ which depends only on $x$, such that $q_T(x)=c_x$, and $\Vert r_nT^nx\Vert^{1/n}\le c_x (\forall n\in\mathbb{N})$, then $T$ is normalike.
	\end{corollary}
	
	With these preparations in hand, we can give examples of forward weighted shift operators which are normalike. Denote $\left\{x\in\mathbb{R}:x\ge0\right\}$ by $\mathbb{R}_+$.

	\begin{lemma}\cite{MR361899}\label{norm}
		Let $T$ be a forward weighted shift operator on $\ell^p$ with weight sequence $\{\omega_k\}_{k=0}^{\infty}\subset\mathbb{R}_+$. Then for every $n\in\mathbb{N}$,
		$$\Vert T^n\Vert=\sup_{k\ge0}\prod\limits_{j=1}^{n}\omega_{k+j-1}.$$
	\end{lemma}

	\begin{proposition}\label{k00}
		Suppose that $T$ is a forward weighted shift operator on $\ell^p$ with weight sequence $\left\{a^{-k}\right\}_{k=0}^{\infty},\ a>1$. Then $\sigma_q(T)=\left\{a^{-k}:k\ge0\right\}$.
	\end{proposition}
	
	\begin{proof}
		Fix $n\in\mathbb{N}$ and $x=\displaystyle\sum\limits_{k=0}^\infty\alpha_k e_k\in\ell^p$ with $\Vert x\Vert=1$. We have $$T^nx=\displaystyle\sum\limits_{k=0}^\infty\alpha_k \left(\prod_{j=1}^{n}\omega_{k+j-1}\right)e_{n+k}=\displaystyle\sum\limits_{k=0}^\infty\alpha_k\cdot a^{-n(2k+n-1)/2}\cdot e_{n+k},$$
		and then$$\Vert T^nx\Vert=\left(\sum_{k=0}^{\infty}|\alpha_k|^p\cdot a^{-pn(2k+n-1)/2}\right)^{1/p}=a^{-n(n-1)/2}\left(\sum_{k=0}^{\infty}|\alpha_k|^p\cdot a^{-pkn}\right)^{1/p}.$$Take $r_n=\Vert T^n\Vert^{-1}$. It is easy to calculate that $\Vert T^n\Vert=a^{-n(n-1)/2}$ by Lemma \ref{norm}, so $$\Vert r_nT^nx\Vert^{1/n}=\left(\sum_{k=0}^{\infty}|\alpha_k|^p\cdot a^{-pkn}\right)^{1/pn}.$$
	
		Set $\Lambda=\left\{k\ge0:\alpha_k\neq0\right\}$, and then $\Lambda\neq\varnothing$ since $\Vert x\Vert=1$. Let $k_0=\min\Lambda$. We have $x=\displaystyle\sum\limits_{k=k_0}^\infty\alpha_k e_k$ and $\alpha_{k_0}\neq0$.
		On one hand,$$\displaystyle\Vert r_nT^nx\Vert^{1/n}=\left(\sum_{k=k_0}^{\infty}|\alpha_k|^p\cdot a^{-pkn}\right)^{1/pn}\le\left(\sum_{k=k_0}^{\infty}|\alpha_k|^p\cdot a^{-pk_0n}\right)^{1/pn}=\Vert x\Vert^{1/n}a^{-k_0}=a^{-k_0}.$$
		On the other hand,$$\displaystyle\Vert r_nT^nx\Vert^{1/n}=\left(\sum_{k=k_0}^{\infty}|\alpha_k|^p\cdot a^{-pkn}\right)^{1/pn}\ge\left(|\alpha_{k_0}|^p\cdot a^{-pk_0n}\right)^{1/pn}=|\alpha_{k_0}|^{1/n}a^{-k_0}.$$
Therefore, $q_T(x)=\displaystyle\limsup_{n\to\infty}\Vert r_nT^nx\Vert^{1/n}=a^{-k_0}$. It follows that $\sigma_q(T)\subset\left\{a^{-k}:k\ge0\right\}$.

At the same time, by the above calculation, we can obtain $q_T(e_k)=a^{-k}$ for every $k\ge0$, so $\sigma_q(T)\supset\left\{a^{-k}:k\ge0\right\}$.

Hence, $\sigma_q(T)=\left\{a^{-k}:k\ge0\right\}$.
	\end{proof}
	
	According to Proposition \ref{k00}, one can see that the quasi-spectrum is not necessarily closed, but it is indeed right-closed.
	
	\begin{proposition}\label{example}
		Let $T$ be as in Proposition \ref{k00}. Then $T$ is normalike.
	\end{proposition}
	
	\begin{proof}
		By Proposition \ref{k00}, for every $x\in X$ with $\|x\|=1$, there exists $k_0\ge0$ such that $q_T(x)=a^{-k_0}$ and $\displaystyle\Vert r_nT^nx\Vert^{1/n}\le a^{-k_0},\ \forall n\in\mathbb{N}$. Then it follows from Corollary \ref{nnnn} that $T$ is normalike.
	\end{proof}

 Recall that a bounded linear operator $T$ on a Hilbert space is said to be hyponormal if $T^\ast T-TT^\ast $ is positive. And if $T$ is hyponormal, then $\|T\|=r(T)$ (see \cite{MR149282}). Now, let $T$ be as in Proposition \ref{k00} and $p=2$. Clearly, $T$ is quasinilpotent, and thus $T$ is not hyponormal. Combining Proposition \ref{k00} and Proposition \ref{example}, one can see that, not only is $T$ a non-hyponormal normalike operator, but also its quasi-spectrum is not a singleton set, even though $T$ is quasinilpotent. This corresponds to facts \ref{sub} and \ref{a-k}.
\vspace{\baselineskip}

Next, we state a difference between the quasi-spectrum and the spectrum. Denote $\left\{\bar{z}:z\in E\right\}$ by $E^\ast$ for $E\subset\mathbb{C}$. Recall that $\sigma(T^\ast)=\sigma(T)^\ast$ for every bounded linear operator $T$ on a complex Hilbert space $\mathcal{H}$. However, this equality does not necessarily hold for the quasi-spectrum. To clarify this, we directly investigate the quasi-spectra for backward weighted shift operators.

\begin{definition}\cite{MR361899}
	Let $\{\omega_k\}_{k=0}^{\infty}\subset\mathbb{R}_+$ be bounded and $1\le p<+\infty$. $T$ is called a backward weighted shift operator on $\ell^p$ if $Te_k=\omega_{k-1}e_{k-1},\ k=1,2,\cdots,\ Te_0=0$.
\end{definition}

\begin{lemma}\cite{MR361899}
	Let $T$ be a backward weighted shift operator on $\ell^p$ with weight sequence $\{\omega_k\}_{k=0}^{\infty}\subset\mathbb{C}$. Then
	$$\Vert T^n\Vert=\sup_{k\ge0}\prod\limits_{j=1}^{n}|\omega_{k+j-1}|,\ \forall n\ge1.$$
\end{lemma}

\begin{lemma}\label{ba}
	Let $T$ be a non-nilpotent backward weighted shift operator on $\ell^p$ with weight sequence $\{\omega_k\}_{k=0}^{\infty}\subset\mathbb{C}$. Then $$q_T(x)=\limsup_{n\to\infty}\left[\dfrac{\left(\sum\limits_{k=0}^\infty|\alpha_{k+n}|^p\prod\limits_{j=1}^{n}|\omega_{k+j-1}|^p\right)^{1/p}}{\sup\limits_{k\ge0}\prod\limits_{j=1}^{n}|\omega_{k+j-1}|}\right]^{1/n},$$where $x=\displaystyle\sum\limits_{k=0}^\infty\alpha_k e_k\in\ell^p$.
\end{lemma}
\begin{proof}
	Suppose $x=\sum_{k=0}^\infty\alpha_k e_k\in\ell^p$. Then$$T^nx=\sum_{k=n}^\infty\left(\alpha_k\prod_{j=1}^{n}|\omega_{k-j}|\right)e_{k-n}=\sum_{k=0}^\infty\left(\alpha_{k+n}\prod_{j=1}^{n}|\omega_{k+j-1}|\right)e_{k},$$
	and hence $$\Vert T^nx\Vert=\left(\sum\limits_{k=0}^\infty|\alpha_{k+n}|^p\prod\limits_{j=1}^{n}|\omega_{k+j-1}|^p\right)^{1/p}.$$
	
	Therefore, $$q_T(x)=\limsup_{n\to\infty}\left[\dfrac{\left(\sum\limits_{k=0}^\infty|\alpha_{k+n}|^p\prod\limits_{j=1}^{n}|\omega_{k+j-1}|^p\right)^{1/p}}{\sup\limits_{k\ge0}\prod\limits_{j=1}^{n}|\omega_{k+j-1}|}\right]^{1/n}.$$

\end{proof}

\begin{corollary}\label{0in}
	Let $T$ be a non-nilpotent backward weighted shift operator on $\ell^p$. Then $q_T(e_k)=0$ for all $k$. Hence, $0\in\sigma_q(T)$.
\end{corollary}

\begin{proposition}\label{back}
	Let $T$ be a backward weighted shift operator on $\ell^p$ with weights 
	$\{\omega_k\}_{k=0}^{\infty}\subset\mathbb{R}_+$, and suppose that 
	$\{\omega_k\}_{k=0}^{\infty}$ decreases monotonically to $0$. Then 
	$\sigma_q(T)=[0,1]$.
\end{proposition}

\begin{proof}
	Let
	\[
	x_r=\sum_{k=0}^{\infty}\frac{e_k}{r^k},\quad r>1.
	\]That is $\alpha_k=r^{-k}$.
	Then for every $n\ge 1$,
	\[
	\left(\sum_{k=0}^{\infty}|\alpha_{k+n}|^p
	\prod_{j=1}^{n}\omega_{k+j-1}^p\right)^{1/p}
	=
	\left(\sum_{k=0}^{\infty}\left(\frac{1}{r^{k+n}}\right)^p
	\prod_{j=1}^{n}\omega_{k+j-1}^p\right)^{1/p}
	=
	\frac{1}{r^n}
	\left(\sum_{k=0}^{\infty}\frac{1}{r^{kp}}
	\prod_{j=1}^{n}\omega_{k+j-1}^p\right)^{1/p},
	\]
	and since $\{\omega_k\}_{k=0}^{\infty}$ decreases monotonically to $0$,
	\[
	\sup_{k\ge 0}\prod_{j=1}^{n}\omega_{k+j-1}
	=
	\prod_{j=1}^{n}\omega_{j-1}.
	\]
	So $$q_T(x_r)=\displaystyle\dfrac{1}{r}\cdot\displaystyle\limsup_{n\to\infty}\left[\displaystyle\dfrac{\displaystyle\left(\sum_{k=0}^{\infty}\frac{1}{r^{kp}}
		\displaystyle\prod_{j=1}^{n}\omega_{k+j-1}^p\right)^{1/p}}{\displaystyle\prod_{j=1}^{n}\omega_{j-1}}\right]^{1/n}$$ by Lemma \ref{ba}.
	
	On the one hand, since $\displaystyle\left(\sum_{k=0}^{\infty}\frac{1}{r^{kp}}
	\displaystyle\prod_{j=1}^{n}\omega_{k+j-1}^p\right)^{1/p}\ge\prod_{j=1}^{n}\omega_{j-1}$,
	\[
	q_T(x_r)
	\ge\dfrac{1}{r}\cdot
	\limsup_{n\to\infty}
	\left(\prod_{j=1}^{n}\omega_{j-1}\right)^{1/n}
	\Big/
	\left(\prod_{j=1}^{n}\omega_{j-1}\right)^{1/n}
	=
	\frac{1}{r}.
	\]
	On the other hand, note that $\displaystyle\left(\sum_{k=0}^{\infty}\frac{1}{r^{kp}}
	\displaystyle\prod_{j=1}^{n}\omega_{k+j-1}^p\right)^{1/p}\le\sum_{k=0}^{\infty}\frac{1}{r^k}\cdot
	\prod_{j=1}^{n}\omega_{k+j-1}$. So
	\[
	q_T(x_r)
	\le\dfrac{1}{r}\cdot
	\limsup_{n\to\infty}
	\left(
	\sum_{k=0}^{\infty}\frac{1}{r^k}\cdot
	\displaystyle\frac{\displaystyle\prod_{j=1}^{n}\omega_{k+j-1}}
	{\displaystyle\prod_{j=1}^{n}\omega_{j-1}}
	\right)^{1/n}
	\le\displaystyle
	\frac{1}{r}\cdot
	\limsup_{n\to\infty}
	\left(
	\sum_{k=0}^{\infty}\frac{1}{r^k}
	\right)^{1/n}
	=
	\frac{1}{r}.
	\]
	Hence \(q_T(x_r)=\frac{1}{r}\). This shows that 
	\(\sigma_q(T)\supset(0,1]\). By Corollary \ref{0in}, \(\sigma_q(T)=[0,1]\).

\end{proof}

 It follows from Proposition \ref{back} that for the operator $T$ in Proposition \ref{k00}, $\sigma_q(T^\ast)\neq\sigma_q(T)^\ast$ when $p=2$. So far, facts \ref{sub} through \ref{spe} have been addressed. We now turn to fact \ref{re}. Before this, we need some preparations.

Ji and Liu in \cite{MR4260659} proved that if a quasinilpotent forward weighted shift $T$ on a Hilbert space is \textit{strongly strictly cyclic} (the definition can be seen in \cite{MR361899}), then the power set of $T$ is $\{1\}$. For the quasi-spectrum, we have the following Theorem.
\begin{theorem}\label{23}
	Suppose that $T$ is a non-nilpotent quasinilpotent forward weighted shift on $\ell^2$. If $T$ is strongly strictly cyclic, then$$\sigma_q(T)=\left\{q_T(e_k):k\ge0\right\}.$$
\end{theorem}

\begin{proof}
	It is pointed out in the proof of Theorem 3.11 in \cite{MR4260659} that for a nonzero vector $x\in\ell^2$, there exist $s_1,\ s_2>0$, and $n_0\in\mathbb{N}$, such that\begin{equation}\label{e1}
		s_1\Vert Ae_{n_0}\Vert\le\Vert Ax\Vert\le s_2\Vert x\Vert\Vert Ae_{n_0}\Vert,\ \forall A\in\mathcal{A}(T|_M),
	\end{equation}
	where $M=\mathcal{H}_{n_0}=\overline{\mathcal{A}(T)x}\in\text{Lat}\ T$. Both $x\in M$ and $e_{n_0}\in M$ are strictly cyclic vector of $T|_M$.
	
	On this basis, fix $n\ge1$, take $A=T^n|_M$, and substitute it into \eqref{e1}. We obtain $$s_1\Vert T^ne_{n_0}\Vert\le\Vert T^nx\Vert\le s_2\Vert x\Vert\Vert T^ne_{n_0}\Vert,$$
	and thus
	\begin{equation}\label{e2}
		s^{1/n}_1\Vert r_nT^ne_{n_0}\Vert^{1/n}\le\Vert r_nT^nx\Vert^{1/n}\le s^{1/n}_2\Vert x\Vert^{1/n}\Vert r_nT^ne_{n_0}\Vert^{1/n},
	\end{equation}
	where $\{r_n\}_{n\in\mathbb{N}}$ is a sequence of positive numbers and $\lim\limits_{n\to\infty}\Vert r_nT^n\Vert^{1/n}=1$. Taking the upper limit as $n\to\infty$ in \eqref{e2} gives $q_T(x)=q_T(e_{n_0})$, and hence $\sigma_q(T)\subset\left\{q_T(e_k):k\ge0\right\}$.
	
	Clearly, $\sigma_q(T)\supset\left\{q_T(e_k):k\ge0\right\}$. So $\sigma_q(T)=\left\{q_T(e_k):k\ge0\right\}$.
\end{proof}

\begin{lemma}\cite{MR361899}\label{ssc}
	Suppose that $T$ is a forward weighted shift with decreasing weight sequence in $\ell^p$ for some $p<+\infty$. Then $T$ is strongly strictly cyclic.
\end{lemma}
\begin{proposition}\label{A}
	Suppose that $A$ is a forward weighted shift operator on $\ell^2$ with weight sequence $\left\{\frac{1}{k+1}\right\}_{k=0}^{\infty}$. Then $\sigma_q(T)=\{1\}$.
\end{proposition}

\begin{proof}
	Since $\{\omega_{k}\}_{k=0}^{\infty}\triangleq\left\{\frac{1}{k+1}\right\}_{k=0}^{\infty}\in\ell^2$, by Lemma \ref{ssc}, $A$ is strongly strictly cyclic. By Theorem \ref{23}, we only need to verify that $q_A(e_k)=1$ for all $k$. Given $k\ge0$ and $n\in\mathbb{N}$, we have$$A^ne_k=\left(\prod\limits_{j=1}^{n}\omega_{k+j-1}\right)e_{n+k},$$ and then$$\Vert A^ne_k\Vert=\prod\limits_{j=1}^{n}\omega_{k+j-1}.$$
	
	Hence,
	$$\begin{aligned}
		q_A(e_k)&=\limsup_{n\to\infty}\left(\dfrac{\prod\limits_{j=1}^{n}\omega_{k+j-1}}{\prod\limits_{j=1}^{n}\omega_{j-1}}\right)^{1/n}\\&=\limsup_{n\to\infty}\left(\dfrac{\dfrac{k!}{(k+n)!}}{\dfrac{1}{n!}}\right)^{1/n}\\
		&=\limsup_{n\to\infty}\left(\dfrac{k!n!}{(k+n)!}\right)^{1/n}=1.
	\end{aligned}$$
	
	So $q_A(e_k)=1,\ k=0,1,2,\cdots$.
\end{proof}

\begin{definition}\cite{MR4654136}
	Given a nonempty subset $E$ of $\mathbb{R}$, define $$\left\lceil E\right\rceil=\left\{\sup F:\emptyset\neq F\subset E,\ F\text{ is bounded}\right\},$$
	and call it the right closure of $E$.
\end{definition}

So $\left\lceil E\right\rceil$ is the smallest right-closed set containing $E$.

\begin{proposition}\label{zk}
	Let $A$ be as in Proposition \ref{A}, and let $\{z_k\}_{k\in\mathbb{N}}\subset[0,1]$ with $z_1=1$. Suppose that $T=\oplus_{k=1}^\infty z_kA$ which is acting on $\mathcal{H}=\oplus_{k=1}^\infty\ell^2$. Then $\sigma_q(T)=\left\lceil\{z_k:k\in\mathbb{N}\}\right\rceil$.
\end{proposition}

\begin{proof}
	It follows from the proof of Theorem B in \cite{MR4654136} that $\Vert T^n\Vert=\Vert A^n\Vert,\ \forall n\in\mathbb{N}$, and that $T$ is quasinilpotent. The aim is to prove $\sigma_q(T)=\sigma\triangleq\left\lceil\{z_k:k\in\mathbb{N}\}\right\rceil$.
	
	\textbullet Fix $m\in\mathbb{N}$. Choose $x=\oplus_{k=1}^\infty x_k\in\mathcal{H}$ such that $x_k\neq0$ if and only if $k=m$. It turns out that $q_T(x)=z_m$. Hence, $\sigma_q(T)\supset\{z_k\}_{k=1}^\infty$. Taking the right closure gives $\sigma_q(T)\supset\sigma$ since $\sigma_q(T)$ is right-closed by Main Theorem \ref{1111}.
	
	\textbullet Fix $x=\oplus_{k=1}^\infty x_k\in\mathcal{H}$ with $\Vert x\Vert=1$. Let $\Lambda=\left\{k\in\mathbb{N}:x_k\neq0\right\}$, and set $\ell=\sup_{k\in\Lambda}\left\{z_k\right\}$, so $\ell\in\sigma$. Then
		$$\begin{aligned}
			q_T(x)&=\limsup_{n\to\infty}\left(\dfrac{\Vert T^nx\Vert}{\Vert T^n\Vert}\right)^{1/n}\\
			&=\limsup_{n\to\infty}\Vert A^n\Vert^{-1/n}\left[\displaystyle\sum_{k\in\Lambda}\Vert (z_kA)^nx_k\Vert^2\right]^{1/2n}\\
			&=\limsup_{n\to\infty}\Vert A^n\Vert^{-1/n}\left(\displaystyle\sum_{k\in\Lambda}z^{2n}_k\Vert A^nx_k\Vert^2\right)^{1/2n}.
		\end{aligned}$$
		
		On the one hand,
		$$\begin{aligned}
			q_T(x)&\le\limsup_{n\to\infty}\Vert A^n\Vert^{-1/n}\left(\displaystyle\sum_{k\in\Lambda}z^{2n}_k\Vert A^n\Vert^2\cdot\Vert x_k\Vert^2\right)^{1/2n}\\
			&=\limsup_{n\to\infty}\left(\displaystyle\sum_{k\in\Lambda}z^{2n}_k\Vert x_k\Vert^2\right)^{1/2n}\\
			&\le\limsup_{n\to\infty}\left(\displaystyle\ell^{2n}\sum_{k\in\Lambda}\Vert x_k\Vert^2\right)^{1/2n}\\
			&\le\limsup_{n\to\infty}\left(\displaystyle\ell^{2n}\right)^{1/2n}=\ell.
		\end{aligned}$$
		
		On the other hand, for every $k_0\in\Lambda$,
		$$\begin{aligned}
		q_T(x)&\ge\limsup_{n\to\infty}\Vert A^n\Vert^{-1/n}\left(\displaystyle z^{2n}_{k_0}\Vert A^nx_{k_0}\Vert^2\right)^{1/2n}\\
		&=z_{k_0}\cdot\limsup_{n\to\infty}\left(\dfrac{\displaystyle \Vert A^nx_{k_0}\Vert}{\Vert A^n\Vert}\right)^{1/n}\\
		&=z_{k_0}\cdot q_A(x_{k_0})\\
		&=z_{k_0}.
		\end{aligned}$$It follows that $q_T(x)\ge\ell$. Hence, $q_T(x)=\ell\in\sigma$. Therefore, $$\sigma_q(T)=\left\{q_T(x):x\in X,\ \Vert x\Vert=1\right\}\subset\sigma.$$

	In summary, $\sigma_q(T)=\sigma$.
	
\end{proof}

	If $\sigma$ is a nonempty right-closed subset of $\mathbb{R}$, there exists a sequence $\{z_n\}_{z\in\mathbb{N}}$ such that its right closure is $\sigma$ (see \cite{MR4654136}). Then by Proposition \ref{zk}, given a right-closed subset $\sigma$ of $[0,1]$ containing 1, there exists a quasinilpotent operator $T\in\mathcal{B}(\mathcal{H})$ such that $\sigma_q(T)=\sigma$.

	\section{The invariant subspace of normalike operators}

In this section, we prove Main Theorem \ref{inv} and Main Theorem \ref{hyp}. Let's see some properties of $q_T(x)$ required for the proof. The first one is very useful. In \cite{MR897643}, it is proved that $r_T(x+y)\le\max\left\{r_T(x),r_T(y)\right\}$ holds for all $x,\ y\in X$. That also holds for $q_T(x)$.

	\begin{lemma}\label{9}
	Let $T\in\mathcal{B}(X)$ be non-nilpotent, let $m\in\mathbb{N}$, and let $x_1,\ x_2,\ \cdots,x_m\in X$. Then $q_T(x)\le\max\limits_{j=1,2,\cdots,m}\{q_T(x_j)\}$, where $x=\sum\limits_{j=1}^{m}x_j$.
\end{lemma}

\begin{proof}
	It suffices to prove the Lemma for $m=2$. Take $x_1,\ x_2\in X$. By the definition of the limit superior, given $\varepsilon>0$, there exists $N>0$ such that for $n>N$, we have
	$$\Vert r_nT^nx_1\Vert^{1/n}\le q_T(x_1)+\varepsilon,\quad \Vert r_nT^nx_2\Vert^{1/n}\le q_T(x_2)+\varepsilon.$$
	Hence,
	$$\begin{aligned}
		q_T(x_1+x_2)&=\limsup_{n\to\infty}\Vert r_nT^n(x_1+x_2)\Vert^{1/n}\\
		&=\limsup_{n\to\infty}\Vert r_nT^nx_1+r_nT^nx_2\Vert^{1/n}\\
		&\le\limsup_{n\to\infty}\left(\Vert r_nT^nx_1\Vert+\Vert r_nT^nx_2\Vert\right)^{1/n}\\
		&\le\limsup_{n\to\infty}\left[\left(q_T(x_1)+\varepsilon\right)^n+\left(q_T(x_2)+\varepsilon\right)^n\right]^{1/n}\\
		&\le\limsup_{n\to\infty}\left[2\left(\max\left\{q_T(x_1),q_T(x_2)\right\}+\varepsilon\right)^n\right]^{1/n}\\
		&=\max\left\{q_T(x_1),q_T(x_2)\right\}+\varepsilon.
	\end{aligned}$$
	
	Since $\varepsilon$ is arbitrary, we obtain $q_T(x_1+x_2)\le\max\left\{q_T(x_1),q_T(x_2)\right\}$.
\end{proof}

\begin{corollary}\label{mani}
	Let $T\in\mathcal{B}(X)$ be non-nilpotent, and let $\tau\in[0,1]$. Then $\left\{x\in X:q_T(x)\le\tau\right\}$ is a linear manifold.
\end{corollary}

\begin{proof}
	In Introduction, we mentioned that $q_T(x)$ is independent of $\|x\|$. Then, combining this property with Lemma \ref{9}, the proof follows immediately.
\end{proof}

In \cite{MR4654136}, it is proved that $k_T(Sx)\le k_T(x)$ holds for $S\in\mathcal{B}(\mathcal{H})$ commuting with $T$.  $q_T(x)$ has the property as well.

\begin{lemma}\label{2}
	Let $T,\ S\in\mathcal{B}(X)$ and $T$ is not nilpotent. If $T$ and $S$ commute, i.e., $ST=TS$, then
	$$q_T(Sx)\le q_T(x),\ \forall x\in X.$$
\end{lemma}

\begin{proof}
	Let $T,\ S\in\mathcal{B}(X)$ such that $ST=TS$. Then for any $x\in X$, we have
	$$\begin{aligned}
		q_T(Sx)&=\limsup_{n\to\infty}\Vert r_nT^nSx\Vert^{1/n}\\
		&=\limsup_{n\to\infty}\Vert r_nST^nx\Vert^{1/n}\\
		&\le\limsup_{n\to\infty}\Vert S\Vert\cdot\Vert r_nT^nx\Vert^{1/n}\\
		&=\limsup_{n\to\infty}\Vert r_nT^nx\Vert^{1/n}\\
		&=q_T(x).
	\end{aligned}$$
\end{proof}

\begin{remark}\label{22}
	By Lemma \ref{2}, given $x\in X$, $q_T(T^kx)\le q_T(x)$ holds for all $k\in\mathbb{N}$.
\end{remark}

\begin{remark}\label{ram}
	For $\tau\in[0,1]$, Corollary \ref{mani} and Lemma \ref{2} show that $\left\{x\in X:q_T(x)\le\tau\right\}$ is a hyperinvariant linear manifold.
\end{remark}

For normalike operators, the following Fatou-type inequality holds.
\begin{lemma}\label{fatou}
	Suppose that $T\in\mathcal{B}(X)$ is normalike. Let $y\in X$, and let $\{y_k\}_{k=1}^\infty\subset X$. If $y_k\to y\ (k\to\infty)$, then $q_T(y)\le\liminf\limits_{k\to\infty}q_T(y_k)$.
\end{lemma}

\begin{proof}
	Let $\{r_n\}_{n\in\mathbb{N}}$ be a sequence of positive numbers with $\lim\limits_{n\to\infty}\|r_nT^n\|^{1/n}=1$. If $q_T(y)=0$, the conclusion is obvious. If $q_T(y)>0$, we aim to prove that $\liminf\limits_{k\to\infty}q_T(y_k)\ge t$ for every $t<q_T(y)$. Denote $\left\{x\in X:\Vert x\Vert=1\right\}$ by $X_1$. Note that
	\[
	\begin{aligned}
		&\left\{x\in X_1:q_T(x)>t\right\}\\
		=&\left\{x\in X_1:\sup_{n\in\mathbb{N}}\|r_nT^n x\|^{1/n}>t\right\}\\
		=&\bigcup\limits_{n=1}^\infty\left\{x\in X_1:\|r_nT^n x\|^{1/n}>t\right\}
	\end{aligned}
	\]
	holds for all $t\in(0,1)$, and the map $x\mapsto\|r_nT^n x\|^{1/n}$ is continuous on $X_1$, so $\left\{x\in X_1:q_T(x)>t\right\}$ is open in $X_1$.
	
	Fix $0<t<q_T(y)$. We have $y/\|y\|\in\left\{x\in X_1:q_T(x)>t\right\}$, since $y\in\left\{x\in X_1:q_T(x)>t\right\}$. Then there exists a neighborhood $U\subset X_1$ of $y/\|y\|$ such that $U\subset\left\{x\in X_1:q_T(x)>t\right\}$. Because $y_k\to y\ (k\to\infty)$, $y_k/\|y_k\|\to y/\|y\|\ (k\to\infty)$. Hence, there exists $K>0$ such that $y_k/\|y_k\|\in U$ for all $k>K$, and thus $q_T(y_k)=q_T(y_k/\|y_k\|)>t$ for all $k>K$. Therefore, $\liminf\limits_{k\to\infty}q_T(y_k)\ge t$. By the arbitrariness of $t$, we obtain $\liminf\limits_{k\to\infty}q_T(y_k)\ge q_T(y)$.
\end{proof}

In Introduction, we mentioned the question: when does \(\sup_{y\in M_x}\{r_T(y)\}=r_T(x)\) hold for every nonzero \(x\in X\)? Then for our version, we may ask: when $T$ is normalike, does \(\sup_{y\in M_x}\{q_T(y)\}=q_T(x)\) hold for every nonzero \(x\in X\)? In fact, it indeed does.
 
\begin{theorem}\label{Mx}
	Suppose that $T\in\mathcal{B}(X)$ is normalike. Let $x\in X\setminus\{0\}$. Then $q_T(y)\le q_T(x)$ for every $y\in M_x$.
\end{theorem}

\begin{proof}
	Fix $x\in X\setminus\{0\}$. Since $q_T(T^mx)\le q_T(x),\ \forall m\ge1$ (by Remark \ref{22}), it follows from Lemma \ref{9} that $q_T(y)\le q_T(x)$ for every $y\in\operatorname{span}_{m\ge0}\{T^mx\}$. For any $y\in M_x$, take a sequence $\{y_k\}_{k=1}^\infty\subset\operatorname{span}_{m\ge0}\{T^mx\}$ such that $y_k\to y(k\to\infty)$. Then we have $q_T(y_k)\le q_T(x)$ for all $k\ge1$. By Lemma \ref{fatou}, we have $\liminf\limits_{k\to\infty}q_T(y_k)\ge q_T(y)$. Consequently, $q_T(x)\ge q_T(y)$.
\end{proof}

Now, we can prove Main Theorem \ref{inv}.

\begin{proof}[\textbf{Proof of Main Theorem \ref{inv}}]
	Take $x,\ y\in X\setminus\{0\}$ with $q_T(x)<1$ and $q_T(y)=1$. Note that $x\in M_x$ and $y\in M_y$. It follows that $M_x\neq\{0\}$ and $M_y\neq\{0\}$. By Theorem \ref{Mx}, we obtain $q_T(z)\le q_T(x)$ for every $z\in M_x$. Hence, $y\notin M_x$, and thus $M_x\neq M_y$. Therefore, at least one of them is not the whole space $X$, meaning that $T$ possesses a nontrivial invariant subspace.
\end{proof}

Furthermore, a nontrivial hyperinvariant subspace can be also constructed. For $\tau\in[0,1]$, denote $\left\{x\in X:q_T(x)\le\tau\right\}$ by $Q_\tau$. Then $Q_\tau$ is a hyperinvariant linear manifold as noted in Remark \ref{ram}.

\begin{proof}[\textbf{Proof of Main Theorem \ref{hyp}}]
	 Our aim is to prove that $Q_{\tau_0}$ is the desired nontrivial hyperinvariant subspace for some $\tau_0\in[0,1)$. So next, we will show $Q_\tau$ is closed for every $\tau\in[0,1]$ and there is a $\tau_0\in[0,1)$ such that $Q_{\tau_0}$ is neither equal to $\{0\}$ nor $X$.
	
	Let $\{y_k\}_{k\in\mathbb{N}}$ be an arbitrary convergent sequence consisting of pairwise distinct elements in $Q_\tau$. Then $q_T(y_k)\le\tau$ for all $k$. Suppose that $\{y_k\}_{k\in\mathbb{N}}$ converges to $y\in X$. By Lemma \ref{fatou}, $q_T(y)\le\liminf\limits_{k\to\infty}q_T(y_k)\le\tau$, which implies $y\in Q_\tau$. Thus $Q_\tau$ is closed.
	
	At last, choose any $\tau_0\in \sigma_q(T)\setminus\{1\}$, and then there is $x_0\in X\setminus\{0\}$ such that $q_T(x_0)=\tau_0$, which means $Q_{\tau_0}\neq\{0\}$. Also, $Q_{\tau_0}\neq X$ by Main Theorem \ref{1111}, obviously.
	
	Therefore, $Q_{\tau_0}$ is a nontrivial hyperinvariant subspace, and we complete the proof.
\end{proof}

\section*{Conflict of Interest}
On behalf of all authors, the corresponding author states that there is no conflict of interest.

	\bibliographystyle{plainurl}
	
	\bibliography{reference}
\end{document}